\documentclass[12pt]{article}
\usepackage{stmaryrd}
\usepackage{mathrsfs}
\usepackage{amsmath}
\usepackage{amsmath,amsthm,amssymb}
\usepackage{amsfonts}
\usepackage{latexsym}
\date{}

\begin{document}
\title{Tilting mutation for $m$-replicated algebras$^\star$}
\author{{\small  Hongbo Lv, Shunhua Zhang$^*$}\\
{\small  Department of Mathematics,\ Shandong University,\ Jinan
250100, P. R. China }}

\pagenumbering{arabic}

\maketitle
\begin{center}
 \begin{minipage}{120mm}
   \small\rm
   {\bf  Abstract}\ \ Let $A$ be a finite dimensional
   hereditary algebra over an algebraically
closed field $k$, $A^{(m)}$ be the $m$-replicated algebra of  $A$
and $\mathscr{C}_{m}(A)$ be the $m$-cluster category of $ A$. We
investigate properties of  complements to a faithful almost complete
tilting $A^{(m)}$-module and prove that the $m$-cluster mutation in
$\mathscr{C}_{m}(A)$ can be realized in ${\rm mod}\ A^{(m)}$, which
generalizes corresponding results on duplicated algebras established
in [Z1].

\end{minipage}
\end{center}
\begin{center}
  \begin{minipage}{120mm}
\end{minipage}
\end{center}

\footnote{ $^\star$Supported by the NSF of China (Grant No.
10771112) and of Shandong Province (Grant No. Y2008A05)}

\footnote{ $^*$ Corresponding author}

\footnote{ {\it Email addresses}: lvhongbo356@163.com(H.Lv), \
shzhang@sdu.edu.cn(S.Zhang)}

\section {Introduction}

\vskip 0.2in

Cluster categories were introduced in [BMRRT] and for type $A_{n}$,
also in [CCS], as a categorical model for better understanding of
cluster algebras of Fomin and Zelevinsky in [FZ1, FZ2].  Now,
cluster categories have become a successful model for acyclic
cluster algebras, see the surveys [BM], [Re] for backgrounds and
recent developments of cluster tilting theory. Later, $m$-cluster
categories were introduced in [Th] as a generalization of cluster
categories. Another good interpretation of $m$-cluster category and
its tilting objects is the $m$-replicated algebras, see [ABST2] and
also see [ABST1] for the case of $m=1$.

\vskip 0.2in

Throughout this paper, we always assume that $A$ is a finite
dimensional hereditary algebra over an algebraically closed field
$k$. Furthermore, we assume that $A$ has $n$ simple modules and
$n\geq 3$ provided $A$ is representation finite. Let $ A^{(m)}$ be
the $m$-replicated algebra of $A$. Then ${\rm gl. dim} A^{(m)}=
2m+1$.  $ A^{(1)}$ is called duplicated algebra.

\vskip 0.2in

Cluster category $\mathscr{C}(A)$ is the orbit category $\mathcal
{D}^b(A)/(\tau^{-1}[1])$ of a bounded derived category $\mathcal
{D}^b(A)$ of $A$ which is a triangulated category by [K],
$m$-cluster category $\mathscr{C}_m(A)$ is the orbit category
$\mathcal {D}^b(A)/(\tau^{-1}[m])$ which also is a triangulated
categories by [K].

\vskip 0.2in

It is well known that there is a one-to-one correspondence between
basic tilting $ A^{(m)}$-modules with projective dimension at most m
and basic tilting objects in $m$-cluster category $\mathscr{C}_m(A)$
, see [ABST2] and see [ABST1] for $m=1$. This motivates further
investigates in this kind of algebras. Some interesting results were
proved in [LLZ], [Z1] and [Z2], for example,  cluster mutation can
be realized in duplicated algebra (see [Z1]). A faithful almost
complete tilting $A^{(m)}$-module with projective dimension at most
$m$ has exactly $m+1$ non-isomorphic complements with projective
dimension at most $m$ (see [LLZ]). Furthermore,  any partial tilting
$ A^{(m)}$-module admits a complement and partial tilting
$A^{(m)}$-module is tilting if and only if the number of its
non-isomorphic indecomposable summands equals to the rank of
Grothendieck group of $ A^{(m)}$ (see [Z2]).

\vskip 0.2in

The aim of this paper is to investigate further properties of
complements to a faithful almost complete tilting $A^{(m)}$-module
and to prove that $m$-cluster mutation in $\mathscr{C}_{m}(A)$ can
be realized in ${\rm mod}\ A^{(m)}$. This paper is arranged as the
following. In section 2, we collect necessary definitions and basic
facts needed for our research. In section 3, we prove a structure
theorem for complements to a faithful almost complete tilting
$A^{(m)}$-module (see Theorem 3.4), and also show that $m$-cluster
mutation in $\mathscr{C}_{m}(A)$ can be realized as tilting mutation
in mod $A^{(m)}$ (see Theorem 3.9). In section 4, we prove that
complements to a faithful almost tilting $A^{(m)}$-module with
projective dimension at most $m$ induce an AR-$(m+3)$-angle in
$\mathscr{C}_{m}(A)$ in the sense of [IY] (see Theorem 4.2).

\vskip 0.2in

\section {Preliminaries}

\vskip 0.2in

Let $\Lambda$ be an Artin algebra.  We denote by mod $\Lambda$ the
category of all finitely generated right $\Lambda$-modules. The
derived category of bounded complexes of mod $\Lambda$ is denoted by
$\mathcal{D}^{b}( \Lambda)$ and the shift functor by [1]. For a
$\Lambda$-module $M$, we denote by add $M$ the subcategory of mod
$\Lambda$ whose objects are the direct summands of finite direct
sums of copies of $M$ and by $\Omega_{\Lambda}^{-1}M$ the first
cosyzygy which is the cokernel of an injective envelope
$M\hookrightarrow I.$ The projective dimension of $M$ is denoted by
pd $M$, the global dimension of $\Lambda$ by gl.dim $\Lambda$ and
the Auslander-Reiten translation of $\Lambda$ by $\tau_\Lambda$.

\vskip 0.2in

Let $\mathcal{C}$ be a full subcategory of mod $\Lambda$,
$C_{M}\in\mathcal{C}$ and $\varphi :C_M\longrightarrow M$ with
$M\in$ mod $\Lambda$. The morphism $\varphi$ is a right
$\mathcal{C}$-approximation of $M$ if the induced  morphism ${\rm
Hom}(C,C_{M})\longrightarrow {\rm Hom}(C,M)$ is surjective for any
$C\in\mathcal{C}$. A minimal right $\mathcal{C}$-approximation of
$M$ is a right $\mathcal{C}$-approximation which is also a right
minimal morphism, i.e., its restriction to any nonzero summand is
nonzero. The subcategory $\mathcal{C}$ is called contravariantly
finite if any module $M\in$ mod $\Lambda $ admits a (minimal) right
$\mathcal{C}$-approximation. The notions of (minimal) left
$\mathcal{C}$-approximation and of covariantly finite subcategory
are dually defined. It is well known that add $M$ is both a
contravariantly finite subcategory and a covariantly finite
subcategory. We call a morphism $\psi : X\longrightarrow Y$ in
$\mathcal{C}$ is a sink map of $Y$ if $\psi$ is right minimal and
${\rm Hom} (\mathcal{C}, X)\longrightarrow {\rm Rad }(\mathcal{C},
Y)\longrightarrow 0$ is exact. A source map can be defined dually.

\vskip 0.2in

Let $T$ be a $\Lambda$-module.  $T$ is said to be exceptional if
${\rm Ext}^i_\Lambda(T, T)=0$ for all $i\geq 1$.  An exceptional
module $T$ is called a partial tilting module provided ${\rm pd} \
T<\infty$. A partial tilting module $T$ is called a tilting module
if  there exists an exact sequence
$$0\longrightarrow \Lambda\longrightarrow T_0\longrightarrow T_1
\longrightarrow \cdots\longrightarrow T_d\longrightarrow 0$$ with
each $T_i\in {\rm add}\  T$. A partial tilting module $T$ is called
an almost complete tilting module if there exists an indecomposable
$\Lambda$-module $N$ such that $T\oplus N$ is a tilting module.

\vskip 0.2in

From now on, let $A$ be a finite dimensional hereditary algebra over
an algebraically closed field $k$.  The repetitive algebra $\hat{A}$
of $A$ is the infinite matrix algebra
$$\hat{A} =\begin{pmatrix}
                                \ddots &  & 0 &  &  \\
                                 & A_{i-1} &  &  &  \\
                                & Q_{i} & A_{i} &  &  \\
                                &  & Q_{i+1} & A_{i+1} &  \\
                                & 0 & & \ddots \\
                             \end{pmatrix}
$$
where matrices have only finitely many non-zero coefficients,
$A_{i}=A$ and $Q_{i}= DA$ for all $i\in \mathbb{Z}$, where $D={\rm
Hom}_k(-,k)$ is the dual functor, all the remaining coefficients are
zero and multiplication is induced from the canonical isomorphisms
$A \otimes_A DA\simeq \ _ADA_A\simeq DA\otimes_AA$ and the zero
morphism $DA\otimes_A DA \longrightarrow 0$ (see [HW] and [H]).

\vskip 0.2in

{\bf  Lemma 2.1.}$^{([H])}$ \ {\it The derived category
$\mathcal{D}^{b}(A)$ is equivalent, as a triangulated category, to
the stable module category ${\rm \underline{mod}}\ \hat{A}$.}

\vskip 0.2in

{\bf Lemma 2.2.}\ {\it Let $M$ be an indecomposable $\hat{A}$-module
which is not projective-injective. Then there exists an
indecomposable $A$-module $N$ such that $M\simeq
\Omega^{-l}_{\hat{A}}N$ for some $l\in \mathbb{Z}$. We denote by $l$
the degree of $M$, that is,  ${\rm deg} M=l$.}

\vskip 0.2in

The $m$-replicated algebra $A^{(m)}$ of $A$ is defined as the
quotient of  the repetitive algebra $\hat{A}$, that is,
$$A^{(m)}=\begin{pmatrix}
                  A_{0} &  &  &  & 0 & \\
                  Q_{1} & A_{1} & &  &  &  \\
                   & Q_{2} & A_{2} &  &  &  \\
                 &  &  & \ddots & \ddots &  \\
                  &  0&  &  & Q_{m} & A_{m} \\
                \end{pmatrix}.
$$

\vskip 0.2in

{\bf Remark.} \ $A^{(1)}$ is the duplicated algebra of $A$ (see
[ABST1]).

\vskip 0.2in

Let $\mathscr{C}_m(A)$ be the $m$-cluster category of $A$.  An
object $X$ in $\mathscr{C}_m(A)$ is said to be exceptional if ${\rm
Ext}^i_{\mathscr{C}_m(A)}(X, X)=0$ for all $i$ with $1\leq i\leq m$
and is called an $m$-cluster tilting object if it is exceptional and
maximal respect to this property. The object $X$ is said to be
almost complete tilting if there is an indecomposable object $Y$
such that $X\oplus Y$ is an $m$-cluster tilting object and $Y$ is
called a complement to $X$. It follows from [ZZ] that, for an almost
complete tilting object $\overline{T}$ in $\mathscr{C}_m(A)$, it has
exactly $m+1$ indecomposable non-isomorphic complements
$\overline{X}_0, \overline{X}_1, \cdots, \overline{X}_m$ and there
are $m+1$ connecting triangles:
$$\ \ \ \ \ \overline{X}_i \stackrel{\overline{f}_i}\longrightarrow
\overline{T}_i \stackrel{\overline{g}_i} \longrightarrow
\overline{X}_{i+1}\longrightarrow \overline{X}_i[1],$$ where
$\overline{f}_i$ is the minimal left ${\rm add}\
\overline{T}$-approximation of $\overline{X}_i$ and $\overline{g}_i$
the minimal right ${\rm add}\ \overline{T}$-approximation of
$\overline{X}_{i+1}$, $i=0,1,\cdots, m.$ We call
$\mu_{\overline{X}_i}(\overline{T}\oplus
\overline{X}_i)=\overline{T}\oplus \overline{X}_{i+1}$ {\it an
$m$-cluster mutation} in direction $\overline{X}_i$.

\vskip 0.2in

The following definition is taken from [IY].

 \vskip 0.2in

{\bf Definition 2.3.}\ {\it Let $M$ be an $m$-cluster tilting object
in $\mathscr{C}_m(A)$. The $(m+3)$-angle
$$X_0 \stackrel{a_0} \longrightarrow M_0 \longrightarrow \cdots
\longrightarrow M_m \stackrel{b_{m+1}}\longrightarrow X_0 $$
induced by $(m+1)$ triangles $$X_i \stackrel{a_i}\longrightarrow M_i
\stackrel{b_{i+1}}\longrightarrow X_{i+1}\longrightarrow X_i[1], \ \
\  i=0,1,\cdots, m$$ is called an AR $(m+3)$-angle if the following
conditions are satisfied:

 \vskip 0.1in

{\rm (1)}\ $X_0$ and $M_i\ (0\leq i\leq m)$ all belong to ${\rm
add}\ M$;

 \vskip 0.1in

{\rm (2)}\ $a_0$ is a source map of $X_0$ in ${\rm add\ }M$ and
$b_{m+1}$ is a sink map of $X_{m+1}=X_0$ in ${\rm add\ }M$;

 \vskip 0.1in

{\rm (3)}\ $a_i$ is a minimal left ${\rm add \ }M$-approximation of
$X_i$ for $1\leq i\leq m$;

 \vskip 0.1in

{\rm (4)}\ $b_i$ is a minimal right ${\rm add \ }M$-approximation of
$X_i$ for $1\leq i\leq m$. }

\vskip 0.2in

We denote by $\pi$ the following composition functor,
$$
\pi:\ {\rm
mod}\ A^{(m)} \hookrightarrow {\rm mod }\ \hat{A}\twoheadrightarrow
{\underline{{\rm mod\ }}}\hat{A} \cong
\mathcal{D}^b(A)\twoheadrightarrow \mathscr{C}_m(A).
$$
By abuse of notation, we often denote objects and modules by the
same letter even when they are considered as objects in different
categories.

\vskip 0.2in

We follow the standard terminology and notation used in the
representation theory of algebras, see [ARS],[H] and [Ri].

\vskip 0.2in

\section {Tilting mutation in ${\rm mod}\ A^{(m)}$}

\vskip 0.2in

The following lemmas are useful and can be easily proved.

\vskip 0.2in

{\bf Lemma 3.1.}\ {\it  For $\hat{A}$-modules $X$ and $Y$,
$${\rm Ext}^s_{\hat{A}}(X,Y)\simeq \underline{{\rm
Hom}}_{\hat{A}}(X, \Omega^{-s}_{\hat{A}}Y).$$}

\vskip 0.2in

{\bf Lemma 3.2.}\ {\it Let $M$ be an indecomposable exceptional
$\hat{A}$-module, which is not projective-injective. Then ${\rm
End}_{\hat{A}}(M)=k.$}

\vskip 0.2in

Let $T$ be a faithful almost complete tilting $A^{(m)}$-module.
According to [LLZ,Z2], we know that $T$ has $t+1$ non-isomorphic
indecomposable complements $X_0, \cdots, X_t$ with $2m\leq t\leq
2m+1$ which are connected by $t$ {\it connecting sequences}:
$$
0\longrightarrow X_i\longrightarrow T_i\longrightarrow
X_{i+1}\longrightarrow 0,\ 0\leq i\leq t-1.
$$
It is easy to see that ${\rm Hom}_{A^{(m)}}(X_j,X_i)=0$ and ${\rm
Ext}^1_{A^{(m)}}(X_i,X_j)=0$ provided $j>i$, and that $X_0$ is the
Bongartz-complement to $T$, which means that $X_0$ can not be
generated by any tilting modules $T\oplus X_i$ for $1\leq i\leq
t-1$. For convenience, we also call $X_t$ the sink complement to
$T$, that is, $X_t$ can not be cogenerated by any tilting modules
$T\oplus X_i$ for $1\leq i\leq t-1$.

\vskip 0.2in

{\bf Lemma 3.3.}\ {\it ${\rm End}_{A^{(m)}}(X_i)=k$ for $0\leq i\leq
t.$}

\vskip 0.2in

{\bf Theorem 3.4.}\ {\it  Taking the notation as above. We have that
$$
{\rm dim}_{k}{\rm Ext}^{s}_{A^{(m)}}(X_j,
X_i)=\left\{\begin{array}{cl}
       1,  & \ \ \ {\rm if}\ i+s=j \ {\rm and }\   s\geq 0, \\
       0, & \ \ \ otherwise.
     \end{array}\right.
$$
Furthermore, the $i$-th connecting sequence $0\longrightarrow
X_i\longrightarrow T_i\longrightarrow X_{i+1}\longrightarrow 0$ is a
$k$-basis of ${\rm Ext}^{1}_{A^{(m)}}(X_{i+1}, X_i)$. Moreover, for
any $0\leq i\leq t-1$ and $0\leq i+s\leq t$, the exact sequence
$0\longrightarrow X_i\longrightarrow T_i\longrightarrow T_{i+1}
\longrightarrow\cdots \longrightarrow T_{i+s-1} \longrightarrow
X_{i+s}\longrightarrow 0$ is a $k$-basis of ${\rm
Ext}^s_{A^{(m)}}(X_{i+s}, X_i)$.}

\vskip 0.1in

{\bf Proof.} \ By applying ${\rm Hom}_{A^{(m)}}(-,X_i)$ to the
$j$-th connecting sequence $$0\longrightarrow X_j\longrightarrow
T_j\longrightarrow X_{j+1}\longrightarrow 0,$$ we get ${\rm
Ext}^{j+1-i}_{A^{(m)}}(X_{j+1},X_i)\cong {\rm
Ext}^{j-i}_{A^{(m)}}(X_{j},X_i)$ for $0\leq i< j\leq t-1$.

In particular, we also have an exact sequence
$$
{\rm Hom}_{A^{(m)}}(X_i,X_i)\rightarrow {\rm
Ext}^{1}_{A^{(m)}}(X_{i+1},X_i)\rightarrow {\rm
Ext}^{1}_{A^{(m)}}(T_i,X_i)=0.
$$
Note that ${\rm Ext}^{1}_{A^{(m)}}(X_{i+1},X_i)\neq 0$ and ${\rm
End}_{A^{(m)}}(X_i)=k$  by Lemma 3.3, it follows that ${\rm
Ext}^{1}_{A^{(m)}}(X_{i+1},X_i)\cong {\rm End}_{A^{(m)}}(X_i)=k$,
and that ${\rm Ext}^{j-i}_{A^{(m)}}(X_{j},X_i)=k$ for $0\leq i<j\leq
t.$  We can take the $i$-th connecting sequence$$0\longrightarrow
X_i\longrightarrow T_i\longrightarrow X_{i+1}\longrightarrow 0 $$ as
a $k$-basis of ${\rm Ext}^{1}_{A^{(m)}}(X_{i+1},X_i)$. It is easy to
see that
$$0\longrightarrow X_i\longrightarrow T_i\longrightarrow
T_{i+1}\longrightarrow\cdots \longrightarrow T_{i+s-1}
\longrightarrow X_{i+s}\longrightarrow 0$$ is non-zero in ${\rm
Ext}^{s}_{A^{(m)}}(X_{i+s},X_i)$, which is a $k$-basis of ${\rm
Ext}^{s}_{A^{(m)}}(X_{i+s},X_i)$.

Now we shall show that ${\rm Ext}^{s}_{A^{(m)}}(X_{j},X_i)=0$ for
$s\neq j-i.$

If $s> j-i$,   $ {\rm Ext}^{s}_{A^{(m)}}(X_{j},X_i)\simeq {\rm
Ext}^{s-(j-i)}_{A^{(m)}}(X_{i},X_i)= 0 $ since  $X_i$ is
exceptional.

Now, we claim that ${\rm Hom}_{A^{(m)}}(X_j,X_i)=0$ provided $j>i$.

In fact, if $j>i+1$, ${\rm Hom}_{A^{(m)}}(X_j,X_i)=0$ since ${\rm
deg} X_j > {\rm deg} X_i$.

We only need to prove that ${\rm Hom}_{A^{(m)}}(X_{i+1},X_i)=0$. On
the contrary we assume that ${\rm Hom}_{A^{(m)}}(X_{i+1},X_i)\neq
0$.

Applying ${\rm Hom}_{A^{(m)}}(X_{i+1},-)$ to the $i$-th connecting
sequence
$$
0\rightarrow X_i\rightarrow T_i\rightarrow X_{i+1}\rightarrow 0,
$$
we have an exact sequence
$$
 0 \rightarrow
{\rm Hom}_{A^{(m)}}(X_{i+1}, X_i)\rightarrow {\rm
Hom}_{A^{(m)}}(X_{i+1}, T_i)\rightarrow {\rm
Hom}_{A^{(m)}}(X_{i+1},X_{i+1})\rightarrow {\rm
Ext}^{1}_{A^{(m)}}(X_{i+1},X_i).
$$
It follows that ${\rm
Hom}_{A^{(m)}}(X_{i+1},T_i)\simeq {\rm
Hom}_{A^{(m)}}(X_{i+1},X_i)\neq 0$ since ${\rm
Hom}_{A^{(m)}}(X_{i+1},X_{i+1})\simeq {\rm
Ext}^{1}_{A^{(m)}}(X_{i+1},X_i)\simeq k$.

In particular, the quiver of algebra ${\rm End}_{A^{(m)}}(T\oplus
X_{i+1})$ will have an oriented cycle, which contradicts with that
$T\oplus X_{i+1}$ being a tilting $A^{(m)}$-module. Our claim is
proved.

If $s <j-i$, according to our claim,
$$
{\rm Ext}^{s}_{A^{(m)}}(X_{j},X_i)\simeq {\rm
Ext}^1_{A^{(m)}}(X_{j-(s-1)},X_i)\simeq{\rm Hom}_{A^{(m)}}(X_{j-s},
X_i)=0.
$$
This completes the proof.  \hfill$\Box$

\vskip 0.2in

Now we are going to show that the converse of Theorem 3.4 is partly
true.

\vskip 0.2in

{\bf Lemma 3.5.}  {\it  Let $M$ be an indecomposable
non-injective-projective $\hat{A}$-module which satisfies that ${\rm
Ext}^s_{\hat{A}}(M,M)=0$ for   $s\geq 1$. Then ${\rm
Ext}^s_{\hat{A}}(\Omega^i_{\hat{A}}M,\Omega^i_{\hat{A}}M)=0$, for
any $i\geq 0$.}

\vskip 0.1in

{\bf Proof.} By Lemma 3.1, we have that
$$\begin{array}{lcl}
{\rm Ext}^s_{\hat{A}}(\Omega^i_{\hat{A}}M,\Omega^i_{\hat{A}}M)
&\simeq & \underline{{\rm Hom}}_{\hat{A}}(\Omega^i_{\hat{A}}M,\Omega^{i-s}_{\hat{A}}M) \\
&\simeq & {\rm Hom}_{\mathcal{D}^b(A)}(M[-i], M[s-i]) \\
&\simeq & {\rm Hom}_{\mathcal{D}^b(A)}(M, M[s]) \\
&\simeq & \underline{{\rm Hom}}_{\hat{A}}(M, \Omega^{-s}M)\\
&\simeq & {\rm Ext}^s_{\hat{A}}(M,M)\\
&= & 0.
\end{array}
$$\hfill$\Box$

\vskip 0.2in

{\bf Definition 3.6.}\ {\it A set of indecomposable
non-projective-injective $A^{(m)}$-modules $\{X_0, X_1, \cdots,
X_t\}$ is called a mutation team in {\rm mod} $A^{(m)}$ if it
satisfies Theorem 3.4, i.e.,
$$
{\rm dim}_{k}{\rm Ext}^{s}_{A^{(m)}}(X_j,
X_i)=\left\{\begin{array}{cl}
       1,  & \ \ \ {\rm if}\ i+s=j \ {\rm and }\   s\geq 0, \\
       0, & \ \ \ otherwise.
     \end{array}\right.
$$
and is maximal with respect to this property.}

\vskip 0.2in

{\bf Remark.} \  Every $X_i$ in a mutation team $\{X_0, X_1, \cdots,
X_t\}$ is exceptional with  ${\rm Hom}_{A^{(m)}}(X_j, X_i)=0$ and
${\rm deg} X_j \geq {\rm deg} X_i$ for $j> i\geq 0$, .

\vskip 0.2in

{\bf Lemma 3.7.}\ {\it Let $\{X_0, X_1, \cdots, X_t\}$ be a mutation
team in {\rm mod} $A^{(m)}$. Then

\vskip 0.1in

(1)\ deg $X_0=0$.

\vskip 0.1in

(2)\ For any $0\leq i\leq t-1$, $0\leq {\rm deg}\ X_{i+1}-{\rm deg}\
X_{i}\leq 1$.

\vskip 0.1in

(3)\ There are at most two elements in $\{X_0, X_1, \cdots, X_t\}$
with same degree.

\vskip 0.1in

(4)\ $2m\leq t\leq 2m+1$.}

\vskip 0.1in

{\bf Proof.}\ (1) Suppose  that deg $X_0=r\geq 1.$ Without loss of
generality, suppose that deg $X_0= 1.$ Then there  is a non-split
exact sequence $$0\longrightarrow \Omega_{A^{(m)}}X_0\longrightarrow
I\longrightarrow X_0\longrightarrow 0,$$ where $I\longrightarrow
X_0\longrightarrow 0$ is a projective  cover and $I$ is
projective-injective. Clearly,
$\Omega_{A^{(m)}}X_0=\Omega_{\hat{A}}X_0$. Since $X_0$ is
exceptional, it follows from Lemma 3.5 that $\Omega_{A^{(m)}}X_0$ is
also exceptional. Applying  ${\rm Hom}_{A^{(m)}}(-,
\Omega_{A^{(m)}}X_0)$ to the sequence above, we have
$$
\cdots\longrightarrow {\rm Hom}_{A^{(m)}}(I,
\Omega_{A^{(m)}}X_0)\longrightarrow {\rm
Hom}_{A^{(m)}}(\Omega_{A^{(m)}}X_0, \Omega_{A^{(m)}}X_0)$$
$$\longrightarrow {\rm Ext^1}_{A^{(m)}}(X_0,
\Omega_{A^{(m)}}X_0)\longrightarrow  {\rm Ext^1}_{A^{(m)}}(I,
\Omega_{A^{(m)}}X_0)\longrightarrow \cdots .
$$
It is easy to see that ${\rm Hom}_{A^{(m)}}(I,
\Omega_{A^{(m)}}X_0)=0$ and that ${\rm Ext^1}_{A^{(m)}}(I,
\Omega_{A^{(m)}}X_0)=0$. By Lemma 3.2, ${\rm
Hom}_{A^{(m)}}(\Omega_{A^{(m)}}X_0, \Omega_{A^{(m)}}X_0)\cong {\rm
Hom}_{\hat{A}}(\Omega_{A^{(m)}}X_0, \Omega_{A^{(m)}}X_0)\cong k$,
which implies that ${\rm Ext^1}_{A^{(m)}}(X_0,
\Omega_{A^{(m)}}X_0)\cong k.$

Denote $\Omega_{A^{(m)}}X_0$ by $X_{-1}$. Then we have that
$$
{\rm dim}_{k}{\rm Ext}^{s+1}_{A^{(m)}}(X_j,
X_{-1})=\left\{\begin{array}{cl}
       1  & \ \ {\rm if}\ s=j \ {\rm and }\  s\geq 0,  \\
       0 & \ \ \ otherwise.
     \end{array}\right.
$$
It follows that $\{X_{-1}, X_0, \cdots, X_t\}$ form a mutation team
in ${\rm mod}\ A^{(m)}$, which is a contradiction. This completes
the proof of (1).

\vskip 0.1in

(2)  Suppose that ${\rm deg}\ X_{i}=r$ and ${\rm deg}\ X_{i+1}=r+p$,
that is, there are two indecomposable $A$-modules $M$ and $N$ such
that $X_i\cong \Omega^{-r}_{A^{(m)}}M$ and $X_{i+1}\cong
\Omega^{-(r+p)}_{A^{(m)}}N$. Then  we see that
$$\begin{array}{rrl}
    k & \cong & {\rm Ext}^1_{A^{(m)}}(X_{i+1}, X_i) \\
      & \cong & \underline{{\rm Hom}}_{A^{(m)}}(X_{i+1}, \Omega^{-1}_{A^{(m)}}X_i) \\
      & \cong & {\rm Hom}_{\mathcal{D}^b(A)}(X_{i+1}, X_i[1]) \\
      & \cong & {\rm Hom}_{\mathcal{D}^b(A)}(N, M[1-p]) \\

  \end{array}
$$It follows that $p=0$ or $p=1.$ This finishes the proof of (2).

\vskip 0.1in

(3) For $p\geq 2$ and $0\leq i\leq t$, we claim that ${\rm deg} X_i
\neq {\rm deg} X_{i+p}$.

Otherwise,  ${\rm deg}\ X_i={\rm deg}\ X_{i+p}$ and $p\geq 2$ imply
that ${\rm Ext}_{\mathcal{D}^{b}( A)}^{p}(X_{i+p}, X_i)=0$.

On the other hand,  we have that
$$\begin{array}{rlr}
   k& \cong{\rm
Ext}_{A^{(m)}}^{p}(X_{i+p}, X_i) \\
      & \cong{\underline{\rm
Hom}}_{A^{(m)}}(X_{i+p}, \Omega_{A^{(m)}}^{-p}X_i)&\ \ \ \ \ \  \\
   & \cong{\rm
Hom}_{\mathcal{D}^{b}( A)}(X_{i+p}, X_i[p]) &\ \ \ \ \ \ \\
   & \cong{\rm
Ext}_{\mathcal{D}^{b}( A)}^{p}(X_{i+p}, X_i),
\end{array}$$

which is a contradiction.

Now we suppose that $X_i,\ X_{i+1}$ and $X_j,\ X_{j+1}$ are the
first four elements in $\{X_0, X_1, \cdots, X_t\}$ such that ${\rm
deg}\ X_i={\rm deg}\ X_{i+1}$ and ${\rm deg}\ X_j={\rm deg}\
X_{j+1}$. It is easy to see that $j>i+1$. According to (1) and (2),
we have that ${\rm deg}\ X_i={\rm deg}\ X_{i+1}=i$ and ${\rm deg}\
X_{j+1}={\rm deg}\ X_j=j-1$ and ${\rm deg}\ X_{j+1}={\rm deg}\
\Omega^{j-1}_{A^{(m)}} X_0$.

On the other hand,
$$\begin{array}{rlr}
   k& ={\rm
Ext}_{A^{(m)}}^{j+1}(X_{j+1}, X_0) \\
      & \cong{\underline{\rm
Hom}}_{A^{(m)}}(X_{j+1}, \Omega_{A^{(m)}}^{-(j+1)}X_0)&\ \ \ \ \ \  \\
   & \cong{\rm
Hom}_{\mathcal{D}^{b}( A)}(X_{j+1}, X_0[j+1]) &\ \ \ \ \ \ \\
   & \cong{\rm
Ext}_{\mathcal{D}^{b}( A)}^{2}(X_{j+1}, X_0[j-1]),
\end{array}$$
which is a contradiction. This completes the proof of (3).

\vskip 0.1in

(4)\ Since ${\rm gl. dim} A^{(m)}= 2m+1$, the consequence follows
from (1), (2) and (3).       \hfill$\Box$

\vskip 0.2in

Recall from [ABST2], the $m$-left part $\mathcal {L}_m(A^{(m)})$ of
${\rm mod }\ A^{(m)}$ consists of the indecomposable
$A^{(m)}$-modules all of whose predecessors have projective
dimension at most $m$.

\vskip 0.2in

{\bf Corollary 3.8.}\ {\it Let $\{X_0, X_1, \cdots, X_t\}$ be a
mutation team in {\rm mod} $A^{(m)}$, and $\{X_0,X_1,\cdots, X_l\}=
\{X_0, X_1, \cdots, X_t\}\bigcap \mathcal {L}_m(A^{(m)})$ be the
partial mutation team in the $m$-left part of ${\rm mod }\ A^{(m)}$.
Then $m-1\leq l\leq m.$}

\vskip 0.2in

{\bf Theorem 3.9.}\ {\it Let $\mathcal{N}$  be a partial mutation
team in the $m$-left part of \ ${\rm mod} A^{(m)}$. Assume that
$\mathcal{N}$ has exactly $m+1$ elements $\{X_0, \cdots, X_m\}$.
Then there exists a faithful almost complete tilting
$A^{(m)}$-module $T$ such that $T\oplus X_i, 0\leq i\leq t $ are all
tilting $A^{(m)}$-modules.}

 \vskip 0.1in

{\bf Proof.}\  The case of $m=1$ has been proved in [Z1] and then we
assume that $m\geq 2$. We only need to prove that
$${\rm Ext}^{l}_{A^{(m)}}(X_j, X_i)={\rm
Ext}^{l}_{\mathscr{C}_m(A)}(\pi(X_j),\pi( X_i))$$ for $1\leq l\leq
m$ and $0\leq i\leq j\leq m$, that is, to show that $\pi(X_0),
\pi(X_1), \cdots, \pi(X_m)$ form an exchange team in
$\mathscr{C}_m(A)$ in sense of [ZZ]. Then according to Theorem 5.8
in [ZZ], there exists an almost complete tilting object $\pi(T')$,
where $T'$ is a non-projective-injective exceptional
$A^{(m)}-$module, such that $\pi(T') \oplus \pi(X_i)$, $0\leq i\leq
m$, are all $m-$cluster tilting objects. By Theorem 29 in [ABST2],
$T'$ has projective dimension at most $m$ and $T'\oplus P$ is a
faithful almost complete tilting $A^{(m)}-$module, where $P$ is the
direct sum of all indecomposable projective-injective
$A^{(m)}-$modules. Let $T=T'\oplus P.$ Then $T$ is just what we
want.

Firstly, we assume that $i+l=j.$ Then ${\rm dim}_k{\rm
Ext}^l_{A^{(m)}}(X_{i+l}, X_i)=1$. Let
$$0\longrightarrow X_i\longrightarrow T_i\longrightarrow
T_{i+1}\longrightarrow\cdots \longrightarrow T_{i+l-1}
\longrightarrow X_{i+l}\longrightarrow 0$$ be a $k-$basis of ${\rm
Ext}^l_{A^{(m)}}(X_{i+l}, X_i)$ given by a chain of non-split short
exact sequences:
$$
0\longrightarrow X_{i+s}\longrightarrow T_{i+s}\longrightarrow
X_{i+s+1}\longrightarrow 0,\ 0\leq s\leq l-1.
$$
By [H], each
$$
0\longrightarrow X_{i+s}\longrightarrow T_{i+s}\longrightarrow
X_{i+s+1}\longrightarrow 0
$$
gives rise to a triangle
$$
X_{i+s}\longrightarrow \overline{T}_{i+s}\longrightarrow
X_{i+s+1}\longrightarrow X_{i+s}[1]
$$
in $\mathcal{D}^b(A)$, which is non-zero in
$$
{\rm Hom}_{\mathcal{D}^b(A)}(X_{i+s+1},X_{i+s}[1])={\rm Ext
}^1_{\mathcal{D}^b(A)}(X_{i+s+1},X_{i+s}).
$$
Then the induced map $X_{i+l}\longrightarrow
X_{i+l-1}[1]\longrightarrow\cdots
\longrightarrow X_{i}[l]$ is non-zero in \\
${\rm Hom}_{\mathcal{D}^b(A)}(X_{i+l}, X_{i}[l])$ and thus
$$
\underline{{\rm Hom}}_{\hat{A}}(X_{i+l},
\Omega_{\hat{A}}^{-l}X_i)\simeq {\rm
Hom}_{\mathcal{D}^b(A)}(X_{i+l},X_i[l])\neq 0.
$$
By the assumption $ {\rm Ext}_{\hat{A}}^l(X_{i+l}, X_i)\simeq{\rm
Ext}_{A^{(m)}}^l(X_{i+l},X_i)\simeq k $ and Lemma 3.1, we get that
$$
{\rm Hom}_{\mathcal{D}^b(A)}(X_{i+l}, X_{i}[l])\simeq
\underline{{\rm Hom}}_{\hat{A}}(X, \Omega^{-l}_{\hat{A}}Y) \simeq
{\rm Ext}_{\hat{A}}^l(X_{i+l}, X_i)\simeq k.
$$
Since $m\geq 2$,  we have that $ {\rm
Hom}_{\mathcal{D}^b(A)}(X_{i+l}, \tau X_i[l-m])=0$ and that

${\rm Hom}_{\mathcal{D}^b(A)}(X_{i+l}, \tau^{-1} X_i[l+m])=0.$
Therefore
$$\begin{array}{lrl}
{\rm Ext}_{\mathscr{C}_m(A)}^l(\pi(X_{i+l}),\pi(X_i)) &\simeq &
{\rm Hom}_{\mathscr{C}_m(A)}(\pi(X_{i+l}), \pi(X_i)[l]) \\
&\simeq & {\rm Hom}_{\mathcal{D}^b(A)}(\pi(X_{i+l}), \pi(X_i)[l])\\
& &\oplus  {\rm Hom}_{\mathcal{D}^b(A)}(\pi (X_{i+l}), \tau\pi(X_i)[l-m]) \\
& &\oplus  {\rm Hom}_{\mathcal{D}^b(A)}(\pi(X_{i+l}), \tau^{-1} \pi(X_i)[l+m])\\
 &\simeq &  {\rm Hom}_{\mathcal{D}^b(A)}(X_{i+l}, X_i[l]) \\
& & \oplus {\rm Hom}_{\mathcal{D}^b(A)}(X_{i+l}, \tau X_i[l-m]) \\
& &\oplus  {\rm Hom}_{\mathcal{D}^b(A)}(X_{i+l}, \tau^{-1} X_i[l+m]) \\
&\simeq & {\rm Hom}_{\mathcal{D}^b(A)}(X_{i+l}, X_i[l]) \\
&\simeq &   k.
 \end{array}
$$
This finishes the proof for the case of $i+l=j.$

Now we assume that $i+l\neq j$. Without loss of generality, we
assume that $i+l< j$.  Then it is easy to see that ${\rm deg}
X_j-{\rm deg }X_i\geq 1.$ Since
$$
{\rm Ext}_{A^{(m)}}^l(X_j, X_i)\simeq {\rm Ext}_{\hat{A}}^l(X_j,
X_i)\simeq {\underline{{\rm Hom}}}_{\hat{A}}(X_j,
\Omega^{-l}_{\hat{A}}X_i)=0,
$$
we have that
$$
{\rm Hom}_{\mathcal{D}^b(A)}(\pi (X_j), \pi (X_i)[l])\simeq {\rm
Hom}_{\mathcal{D}^b(A)}(X_j, X_i[l])\simeq {\underline{\rm
Hom}}_{\hat A}( {X_j}, \Omega^{-l}_{\hat{A}}{X_i})=0.
$$
Since $m\geq 2,$ it is easy to see that
$$
{\rm Hom}_{\mathcal{D}^b(A)}(\pi(X_j), \tau \pi(X_i)[l-m])\cong {\rm
Hom }_{\mathcal{D}^b(A)}(X_j, \tau X_i[l-m])=0.
$$
We now claim that ${\rm Hom }_{\mathcal{D}^b(A)}(\pi(X_j), \tau^{-1}
\pi(X_i)[l+m])=0$. Note that $l+i+m-j\geq 1$ since $l\geq 1$ and
$m\geq j.$ Then if $l+i+m-j\geq 2$,  our claim is true since ${\rm
deg} \pi(X_j) \leq {\rm deg }\tau^{-1}\pi(X_i)[j-i]\leq {\rm deg
}\tau^{-1}\pi(X_i)[l+m]-2$. For the case of $l+i+m-j=1$, that is,
$l=1,\ i=0$ and $j=m$, we will show that
$$
{\rm Hom}_{\mathcal{D}^b(A)}(\pi(X_m),
\tau^{-1}\pi(X_0)[m][1])\simeq {\rm Hom}_{\mathcal{D}^b(A)}(X_m,
\tau^{-1}X_0[m][1])=0.
$$
By Lemma 3.4 in [LLZ], ${\rm pd} _{A^{(m)}}X_m=m  $ and thus ${\rm
deg} X_m$ is either $m-1$ or $m$.

If ${\rm deg} X_m=m-1$, our claim holds because ${\rm deg
}\tau^{-1}X_0[m+1]\geq m+1$.

If ${\rm deg} X_m=m$, there exists an indecomposable projective
$A$-module $P$ such that $X_m=P[m]$. Then
$$
{\rm Hom}_{\mathcal{D}^b(A)}(X_m, \tau^{-1}X_0[m][1])\simeq {\rm
Hom}_{\mathcal{D}^b(A)}(P, \tau^{-1}X_0[1])=0.
$$

By the arguments above, we get that, for the case of $l+i< j,$
$$
\begin{array}{lll}
{\rm Ext}_{\mathscr{C}_m(A)}^l(\pi(X_{i+l}),\pi(X_i))
&\simeq &  {\rm Hom}_{\mathscr{C}_m(A)}(\pi(X_{i+l}), \pi(X_i)[l]) \\
&\simeq & {\rm Hom}_{\mathcal{D}^b(A)}(\pi(X_{i+l}), \pi (X_i)[l])\\
&& \oplus {\rm Hom}_{\mathcal{D}^b(A)}(\pi (X_{i+l}), \tau \pi(X_i)[l-m]) \\
&&  \oplus {\rm Hom}_{\mathcal{D}^b(A)}(\pi(X_{i+l}), \tau^{-1} \pi(X_i)[l+m])\\
 &\simeq & {\rm Hom}_{\mathcal{D}^b(A)}(X_{i+l}, X_i[l])\\
& & \oplus {\rm Hom}_{\mathcal{D}^b(A)}(X_{i+l}, \tau X_i[l-m]) \\
 && \oplus {\rm Hom}_{\mathcal{D}^b(A)}(X_{i+l}, \tau^{-1} X_i[l+m]) \\
& =& 0.
\end{array}
$$
This finishes the proof of the theorem.    \hfill$\Box$

\vskip 0.2in

\section {Relation with AR $(m+3)$-angle in $\mathscr{C}_m(A)$ }

\vskip 0.2in

In this section, we shall give a further explanation about the
relationship between the tilting mutation in ${\rm mod }\ A^{(m)} $
and the $m$-cluster mutation in $\mathscr{C}_m(A).$

\vskip 0.2in

Let $T$ be a faithful almost complete tilting $A^{(m)}$-module with
${\rm pd}_{A^{(m)}}T\leq m$. By [LLZ], $T$ has exactly $m+1$
indecomposable non-isomorphic complements $X_0, \cdots, X_m$ with
projective dimensions at most $m$, which are connected by the long
exact sequence:
$$
(\ast)\ \ 0\longrightarrow X_0\longrightarrow
T_0\longrightarrow T_{1}\longrightarrow\cdots \longrightarrow
T_{m-1} \longrightarrow X_{m}\longrightarrow 0,
$$
where $T_i\in{\rm add}\, T $ for all $0\leq i\leq m-1$, $ X_i={\rm
Coker}\ g_{i-1}$ for $1\leq i\leq m$ and each of the induced
monomorphisms $X_i\hookrightarrow T_i$ is a  minimal left ${\rm
add}\ T$-approximation.  It follows from Theorem 29 in [ABST2] that
$\pi(T)$ is an almost complete $m$-cluster tilting object in
$\mathscr{C}_m(A)$, and that $\pi(X_0), \cdots, \pi(X_m)$ are its
$m+1$ indecomposable non-isomorphic complements, which are connected
by the connecting triangles:
$$
\pi(X_i)\overset{\overline{f_i}}{\longrightarrow}
\pi(T_i)\overset{\overline{g_i}}{\longrightarrow}
\pi(X_{i+1})\longrightarrow \pi(X_i)[1],
$$
where $\overline{f_i} $
is the minimal left ${\rm add }\ \pi(T)$-approximation of $\pi(X_i)$
and $\overline{g_i}$ is the minimal right ${\rm add }\
\pi(T)$-approximation of $\pi(X_{i+1})$. Then we have a long angle:
$$
(\ast\ast)\ \  \pi(X_0)\longrightarrow \pi(T_0)\longrightarrow \pi(T_{1})\longrightarrow\cdots
\longrightarrow \pi(T_{m}) \longrightarrow \pi(X_{0}).
$$
It is easy to see that $(\ast\ast)$ is an AR $(m+3)$-angle for
$m\geq 2$ by Corollary 4.4 in [ZZ] and for $m=1$ by Lemma 6.13 in
[BMRRT]. Now, we want to show that $(\ast\ast)$ is induced by
$(\ast)$.

\vskip 0.2in

{\bf Lemma 4.1}\ {\it For all $0\leq i\leq m-1$, the $i$-th
connecting sequence $$0\longrightarrow
X_i\overset{f_i}{\longrightarrow} T_i\overset{g_i}{\longrightarrow}
X_{i+1}\longrightarrow 0$$ in ${\rm mod}\ A^{(m)}$ induces a
triangle in $\mathscr{C}_m(A)$:
$$\pi(X_i)\overset{\overline{f}_i}{\longrightarrow}
\pi(T_i)\overset{\overline{g}_i}{\longrightarrow}
\pi(X_{i+1})\longrightarrow \pi(X_i)[1],$$ where $\overline{f}_i $
is the minimal left ${\rm add }\ \pi(T)$-approximation of $\pi(X_i)$
and  $\overline{g}_i$ is the minimal right ${\rm add }\
\pi(T)$-approximation of $\pi(X_{i+1})$.}

\vskip 0.1in

{\bf Proof.}\ By [H], the short exact sequence in ${\rm mod}\
A^{(m)}$ (also in ${\rm mod}\ \hat{A}$)
$$
 0\longrightarrow
X_i\overset{f_i}{\longrightarrow} T_i\overset{g_i}{\longrightarrow}
X_{i+1}\longrightarrow 0
$$
gives rise to a triangle  $ X_i{\longrightarrow}
\overline{T}_i{\longrightarrow} X_{i+1}\longrightarrow X_i[1]$ in
$\mathcal{D}^b(A)$ and hence a triangle in $\mathscr{C}_m(A)$:
$$
\pi(X_i)\overset{\overline{f}_i}{\longrightarrow}
\pi(T_i)\overset{\overline{g}_i}{\longrightarrow}
\pi(X_{i+1})\longrightarrow \pi(X_i)[1].
$$
Note that $\pi$ is an exact functor and ${\rm deg}X_i\leq m$ for all
$0\leq i\leq m.$ By Lemma 2.1 and Lemma 3.1, non-split sequence
$0\longrightarrow X_i\overset{f_i}{\longrightarrow}
T_i\overset{g_i}{\longrightarrow} X_{i+1}\longrightarrow 0$ induces
a non-split triangle $ X_i{\longrightarrow} T_i{\longrightarrow}
X_{i+1}\longrightarrow X_i[1]$  in $\mathcal{D}^b(A)$. Moreover,
$\pi(X_i)\overset{\overline{f}_i}{\longrightarrow}
\pi(T_i)\overset{\overline{g}_i}{\longrightarrow}
\pi(X_{i+1})\longrightarrow \pi(X_i)[1]$ is not-split in
$\mathscr{C}_m(A)$.

It follows from [ZZ] that ${\rm dim}_k{\rm
Ext}^1_{\mathscr{C}_m(A)}(\pi(X_{i+1}), \pi(X_i))=1$.

Let $\pi(X_i)\overset{\alpha}{\longrightarrow}
E_i\overset{\beta}{\longrightarrow} \pi(X_{i+1})\longrightarrow
\pi(X_i)[1]$ be a basis of ${\rm
Ext}^1_{\mathscr{C}_m(A)}(\pi(X_{i+1}), \pi(X_i))$, where $\alpha $
(resp. $\beta$) is the minimal left (resp.  right) ${\rm
add}\pi(T)-$approximation of $\pi(X_i)$ (resp. $\pi(X_{i+1})$). Then
$ \pi(X_i)\overset{\overline{f}_i}{\longrightarrow}
\pi(T_i)\overset{\overline{g}_i}{\longrightarrow}
\pi(X_{i+1})\longrightarrow \pi(X_i)[1]$ is isomorphic to this basis
and hence $\overline{f}_i $ is the minimal left ${\rm add\
}\pi(T)$-approximation of $\pi(X_i)$ and $\overline{g}_i$ is the
minimal right ${\rm add }\ \pi(T)$-approximation of $\pi(X_{i+1})$.
\hfill$\Box$

\vskip 0.2in

{\bf Theorem 4.2.}\ {\it Let $T$ be a faithful almost complete
tilting $A^{(m)}$-module with projective dimension at most $m$ and
$X_0, \cdots, X_m$ be its indecomposable non-isomorphic complements
with projective dimension at most $m$. Then the induced
$(m+3)$-angle in  $\mathscr{C}_m(A)$ by the $m+1$ connecting
sequences is just the AR ($m+3$)-angle in the sense of [IY].}

\vskip 0.1in

{\bf Proof.}\  By Lemma 4.1, the long exact sequence in {\rm mod}
$A^{(m)}$
$$(\ast)\ \ 0\longrightarrow X_0\longrightarrow T_0\longrightarrow
T_{1}\longrightarrow\cdots \longrightarrow T_{m-1} \longrightarrow
X_{m}\longrightarrow 0$$ induces a long angle in $\mathscr{C}_m(A)$
$$ (1\ast)\ \ \pi(X_0)\longrightarrow \pi(T_0)\longrightarrow
\pi(T_{1})\longrightarrow\cdots \longrightarrow \pi(T_{m-1})
\longrightarrow \pi(X_{m}).$$ Let $\pi(X_{m})\longrightarrow
\pi(T_{m})$ be the minimal left ${\rm add }\ \pi(T)$-approximation,
which induces a triangle in $\mathscr{C}_m(A)$
$$(2\ast)\ \ \pi(X_m){\longrightarrow}
\pi(T_m){\longrightarrow} \pi(X_{0})\longrightarrow \pi(X_m)[1],$$
where $\pi(T_m){\longrightarrow} \pi(X_{0})$ is the minimal right
${\rm add }\ \pi(T)$-approximation by [BMRRT]. Now we get the long
angle $(\ast\ast)$ by connecting  $(1\ast)$ and $(2\ast)$. By
Corollary 4.4 in [ZZ], the proof is finished. \hfill$\Box$

\begin{description}

\item{[ABST1]}\ I.Assem, T.Br$\ddot{\rm u}$stle, R.Schiffer,
G.Todorov,
 Cluster categories and duplicated algebras, J.
Algebra. 305 (2006), 548-561.

\item{[ABST2]}\ I.Assem, T.Br$\ddot{\rm u}$stle, R.Schiffer,
G.Todorov,
 $m$-cluster categories and $m$-replicated algebras,
Journal of pure and applied algebra. 212 (2008), 884-901.

\item{[ARS]}\ M.Auslander, I.Reiten, S.O.Smal$\phi$, \ Representation
Theory of Artin Algebras, Cambridge Univ. Press, 1995.

\item{[BM]} A. Buan, R. Marsh, Cluster-tilting theory. Trends in
representation theory of algebras and related topics, Edited by J.
de la Pe$\widetilde{n}$a and R. Bautista.  Contemporary Math.,
406(2006), 1-30.

\item{[BMRRT]}\  A. Buan, R. Marsh, M. Reineke, I. Reiten, G.Todorov,
Tilting theory and cluster combinatorics, Adv.in Math.,
204(2)(2006), 572-618.

\item{[CCS]}\  P. Caldero, F. Chapoton, R. Schiffler,
                   Quivers with relations arising from
                  clusters($A_n$ case).  Trans. Amer. Math. Soc.,
                   358(3)(2006), 1347-1364.

\item{[FZ1]} \ S. Fomin, A. Zelevinsky,
Cluster algebra I: Foundation.  J. Amer. Math. Soc., 15(2002),
497-529.

\item{[FZ2]} \   S. Fomin, A. Zelevinsky,
             Cluster algebra II: Finite type classification.
                Invent.Math.,  154(2003), 63-121.

\item{[H]}\ D.Happel, Triangulated categories in the representation
theory of finite dimensional algebras, Lecture Notes series 119.
Cambridge Univ. Press, 1988.

\item{[HW]}\ D.Hughes, J.Waschb$\ddot{\rm u}$sch. Trivial extensions
of tilted algebras.  Proc.London. Math.Soc., 46(1983), 347-364.

\item{[IY]} O. Iyama, Y. Yoshino, Mutation in triangulated
categories and rigid Cohen-Macaulay modules.  Invent. Math.,
172(2008), 117-168.

\item{[K]}\ B. Keller, Triangulated orbit categories.  Documenta
Math. 10(2005),551-581.

\item{[LLZ]}\ X.Lei, H.Lv, S.Zhang, Complements to the almost
complete tilting $A^{(m)}$-modules, To appear in  Communications in
Algebra.

\item{[Re]} I. Reiten, Tilting theory and cluster algebras.
Preprint.

\item{[Ri]}\ C.M.Ringel, Tame algebras and integral quadratic
forms, Lecture Notes in Math. 1099. Springer Verlag, 1984.

\item{[Th]}\ H.Thomas, Defining an $m$-cluster category. J.Algebra,
318(2007), 37-46.

\item{[Z1]}\ S.Zhang, Tilting mutation and duplicated algebras. To
appear in Communications in Algebra.

\item{[Z2]}\ S.Zhang, Partial tilting modules over $m$-replicated
algebras.  Preprint, arXiv: math.RT/0810.5190, 2008.

\item{[ZZ]} Y. Zhou, B. Zhu, Cluster combinatorics of $d$-cluster
categories. Preprint, arXiv: math.RT/0712.1381, 2007.

\end{description}

\end{document}